\documentclass[smallextended,natbib,runningheads]{svjour3}

\usepackage{bm}
 
\usepackage{amsmath}
\usepackage{amssymb}
\usepackage{theorem}
\usepackage{xspace}
\usepackage{geometry}
\usepackage{array}
\usepackage{subfigure}
\usepackage{epsfig}
\usepackage{hhline}
\usepackage{bbm}
 
\usepackage[utf8]{inputenc}    
\DeclareMathAccent{\widehat}{\mathord}{largesymbols}{"62}
\DeclareMathAccent{\widetilde}{\mathord}{largesymbols}{"65}
\def\pth#1{\left(#1\right)}
\def\acc#1{\left\{#1\right\}}
\def\cro#1{\left[#1\right]}

\def\eeX{\mathbb{X}}

\def\ebo{\textrm{\mathversion{bold}$\mathbf{\beta^0}$\mathversion{normal}}}
\def\uu{\textrm{\mathversion{bold}$\mathbf{1}$\mathversion{normal}}}

\def\eb{\textrm{\mathversion{bold}$\mathbf{\beta}$\mathversion{normal}}}  
\def\ed{\textrm{\mathversion{bold}$\mathbf{\delta}$\mathversion{normal}}}

\def\eU{\textrm{\mathversion{bold}$\mathbf{\Upsilon}$\mathversion{normal}}} 

\def\eE{I\!\!E}
\def\eP{I\!\!P}
\def\e1{1\!\!1}

\def\XX{\textrm{\mathversion{bold}$\mathbf{X}$\mathversion{normal}}}

\def\xx{\textrm{\mathversion{bold}$\mathbf{x}$\mathversion{normal}}}
\def\hh{ \hspace*{0.5cm}}

\theoremstyle{plain}

\newcommand{\beqn}{\begin{eqnarray*}}
\newcommand{\eeqn}{\end{eqnarray*}}
  
\def\ee1{\textrm{\mathversion{bold}$\mathbf{\varepsilon}$\mathversion{normal}}}

\newcommand{\bfz}{{\bf z}}

\def\eu{\mathbf{{u}}}

\newcommand{\N}{\mathbb{N}}
\newcommand{\R}{\mathbb{R}}
\newcommand{\PP}{\mathbb{P}}

\def\argmin{\mathop{\mathrm{arg\,min}}} 

\begin{document}
\title {{\bf Adaptive group LASSO selection in quantile models}}

\author{GABRIELA CIUPERCA}
\institute{GABRIELA  CIUPERCA \at
Universit\'e de Lyon,
Universit\'e Lyon 1, 
CNRS, UMR 5208, Institut Camille Jordan, 
Bat.  Braconnier,
43, blvd du 11 novembre 1918,
F - 69622 Villeurbanne Cedex, France,\\\email{Gabriela.Ciuperca@univ-lyon1.fr}\\
{\it tel: }33(0)4.26.23.45.57, {\it fax: }33(0)4.72.43.16.87}
\date{Received: date / Accepted: date}
\maketitle

\begin{abstract}
 The paper considers a linear model with  grouped  explanatory variables. If the model errors are not with zero mean and bounded variance or if model contains outliers, then the least squares framework is not appropriate. Thus, the quantile regression is an interesting alternative. In order to automatically select the relevant variable groups, we propose and study here the adaptive group LASSO quantile estimator.  
 We establish the sparsity and asymptotic normality of the proposed estimator in two cases: fixed number and divergent number of variable groups. Numerical study by Monte Carlo simulations confirms the theoretical results and illustrates the performance of the proposed estimator.  
\end{abstract}
 
\section{Introduction}
Classically, for the   regression model, the  errors are assumed to be independent, of mean zero  and  bounded variance. Then, the model is estimated by the least squares (LS) method, eventually with a penalty of  LASSO  type  when automatic detection of significant variables is performed. 
If the assumptions on the first two moments of the  model error are not satisfied, then the LS framework breaks down. In this case, an alternative is to consider the quantile regression with a LASSO  type penalty. This is one of the interests of this paper. The quantile regression is robust and allows relaxation of the two first moment conditions of the model error.\\
Often enough in practice,  for example in the variance analysis case,  are considered the regression linear models with  grouped variables. For models with grouped explanatory variables it is more meaningful to identify relevant variable groups instead of individual variables.  If the errors have  Normal distribution, then for detecting the relevant variable groups, the    F-statistic test is used. If the errors are not  Gaussian and if more the number of groups is large, then the F-statistic  test is inappropriate. From where, another interest of this paper: we  consider the quantile process with LASSO type  penalty in order to automatically detect the irrelevant variable groups.  \\
\hh The automatic selection method of the grouped variables  using the  LASSO penalties was introduced by \cite{Yuan-Lin-06} for  gaussian errors, by proposing  the   LASSO group penalty for the process of the error squares sum. Several recent papers have considered group selection using LASSO type penalties. For fixed parameter space and mean zero, finite second moment i.i.d.  model  errors,   \cite{Nardi-Rinaldo-08} established the model selection consistency and asymptotic normality of nonzero group LASSO estimator. The same estimator is studied by  \cite{Nardi-Rinaldo-08} when number of covariates is larger, for particular case of normal errors. For gaussian errors, \cite{Xu-Ghosh-15} realize a Bayesian variable selection by penalization of the error squares sum with Bayesian group LASSO. For this estimation method, the posterior median estimator satisfies the sparsity property.  The adaptive group LASSO estimators, when  the number $p$ of groups is fixed,   was studied by \cite{Wang-Leng-08}. For high-dimensional model, \cite{Wei-Huang-10} studied the selection and estimation properties of the adaptive group LASSO,  but under assumption that the errors are gaussian. Still for the error squares sum penalized with adaptive LASSO penalty, \cite{Zhang-Xiang-15} consider the case of the number of groups $p_n $ converges to infinity when $n \rightarrow \infty$, for  i.i.d. errors $\varepsilon$ such that   $\eE[\varepsilon]=0$ and $Var[\varepsilon]< \infty$. The consistency and asymptotic normality of the parameter estimator are established.  A paper that doesn't consider the LS penalized process, but a process associated to a twice differentiable convex function,  with LASSO penalty, for the  case $p$ large and small $n$ was considered by  \cite{Wang-You-Lian-15}.  When the number of groups can grow at a certain polynomial rate,   the automatic selection property  of variable  groups  for a LS  process with  SCAD penalty has been proven in \cite{Guo-Zhang-Wang-Wu-15}. Automatic selection of the   relevant  variable groups, when $p  $ converges to infinity, has also considered by \cite{Zou-Zhang-09} penalizing the LS process with adaptive elastic-net penalty.  For a review of group selection methods and several applications of these methods the reader can see \cite{Huang-Breheny-Ma-12}.\\
\hh In this paper we consider the model selection problem  and the  estimation in a linear model with $p$ groups of explanatory variables. We propose and study the asymptotic properties of the adaptive group  LASSO  quantile estimator in two  cases:   $p$ fixed and  $p \rightarrow \infty$ as $n \rightarrow \infty$. This  estimator is the minimizer of the quantile process penalized  by an adaptive group LASSO penalty. The oracle properties, i.e. the  automatic selection of significant variables groups and their asymptotic distribution, are proved.\\
\hh The remainder of the  paper is organized as follows. In Section 2 we present the model and  introduce some notations used throughout in this paper. Oracle properties for the adaptive group   LASSO quantile estimator are proved  for  $p$ fixed in Section 3 and for $p \rightarrow \infty$ as $n \rightarrow \infty$ in Section 4.   Section 5 reports some simulation results which illustrate the method interest. We compare the adaptive group LASSO quantile  estimation performance with the adaptive group LASSO least squares  estimations, proposed by \cite{Zhang-Xiang-15}.   All proofs are given in  Section 6.   
\section{Model and  notations}
In this section, we present the statistical model and we also introduce some notations used throughout in the paper.\\
We begin by introducing some general notations. All vectors and matrices are denoted by bold symbols and all vectors are written as  column vectors. For a vector $\textbf{v}$, we denote by $\textbf{v}^t$ its transposed and by $\|\textbf{v} \|$ its Euclidean norm. Notations $ \overset{\cal L} {\underset{n \rightarrow \infty}{\longrightarrow}}$, $ \overset{\PP} {\underset{n \rightarrow \infty}{\longrightarrow}}$ represent the convergence in distribution and in probability, respectively, as $n \rightarrow \infty$.  For a positive definite matrix $\textbf{M}$, we denote by $\lambda_{\min}(\textbf{M)}$  and $\lambda_{\max}(\textbf{M)}$ its the smallest and largest eigenvalues, respectively. \\
 We will also use   the following notations: if $V_n $ and $U_n$ are random variable sequences, $V_n=o_{\eP}(U_n)$ means that $\lim_{n \rightarrow \infty} \eP[|U_n/V_n| > e]=0$ for any $e>0$, $V_n=O_{\eP}(U_n)$ means that there exists a finite $C>0$ such that $\eP[|U_n/V_n| > C]< e$ for any $n$ and $e$. If $V_n $ and $U_n$ are  deterministic sequences,   $V_n=o (U_n)$ means that the sequence $V_n/U_n \rightarrow 0$ for $n \rightarrow\infty$,  $V_n=O (U_n)$ means that the sequence $V_n/U_n $ is bounded for sufficiently large $n$.\\
   Throughout this paper, $C$ will denote generic constant; not depending on  size $n$ which may take different values in different formula or even in different parts of the same formula. The value of $C$ is not of interest. We will also use the notation $\textbf{0}_k$ for the zero k-vector.\\

We consider the following linear model with $p$ groups of explanatory variables:
\begin{equation}
\label{eq1}
Y_i=\sum^p_{j=1}\XX_{ij}^t \eb_j +\varepsilon_i=\eeX_i^t \eb +\varepsilon_i, \qquad  i =1, \cdots, n,
\end{equation}
with $Y_i, \varepsilon_i$  random variables. For each group  $j =1, \cdots, p$, the vector of the parameters is  $\eb_j \equiv (\beta_{j1}, \cdots , \beta_{j d_j}) \in \R^{d_j}$ and the design for observation  $i$ is  $\XX_{ij}$,  a column vector of size $d_j$. The vector with all coefficients is $\eb\equiv (\eb_1, \cdots, \eb_p) $ and for observation $i$,  the vector with all explanatory variables is   $\eeX_i=(\XX_{i1}, \cdots , \XX_{ip})$.  
Denote by $ \eb^0_j=(\beta^0_{j1}, \cdots ,\beta ^0_{j d_j})$ the true value (unknown) of the parameter $\eb_j$. For observation $i$, we denote by $X_{ij,k}$ the $k$th variable of the  $j$th  group.\\
We emphasize that for the $i$th sample, we observe $(Y_i, \eeX_i)$, $i=1, \cdots, n$. \\ 
The relevant groups of explanatory variables correspond to the nonzero vectors. Without loss of generality, on suppose that the first $p_0$ ($p_0 \leq p$) groups of explanatory variables are relevant: 
 \[
 \| \eb^0_j\| \neq 0, \quad \textrm{for all } j \leq p_0 \qquad \textrm{and } \| \eb^0_j\| =0, \quad \textrm{for all } j > p_0,
 \]
where  $\| .\|$ is the Euclidean norm. Let $r$ be the total number of explanatory variables, so $r=\sum^{p}_{j=1}d_j$. We denote by $r^0=\sum^{p_0}_{j=1}d_j$. So, $p_0$ is the number of nonzero true parameter vectors and $r^0$ is the total number of parameters in these nonzero true vectors. \\
 The multi-factor ANOVA model is an example of this model.\\

We introduce now the quantile framework.  For a fixed quantile index $\tau \in (0,1)$, the check function  $\rho_\tau(.): \R \rightarrow \R_{+}$ is  defined by  $\rho_\tau(u)= u(\tau -\e1_{u< 0})$.\\
  The quantile estimator of   $\eb$, is the  minimizer of  the    quantile process associated to model (\ref{eq1}):
  \begin{equation}
  \label{eq2}
  \widetilde{\eb}_n \equiv \argmin_{\eb \in \R^r} \sum^n_{i=1}\rho_\tau(Y_i -\eeX_i^t \eb).
  \end{equation}
  For the particular case $\tau=1/2$ we obtain the median regression and (\ref{eq2}) becomes the least absolute deviations estimator. A great advantage of the quantile framework is that, compared to classical estimation methods that are sensitive to outliers, the quantile method provides more robust estimators. Moreover, the required assumptions to the error moments are relaxed. \\
\hh  The  estimator $\widetilde{\eb}_n=(\widetilde{\eb}_{n;1}, \widetilde{\eb}_{n;2}, \cdots , \widetilde{\eb}_{n;p})$ has as   $d_j$-subvector $\widetilde{\eb}_{n;j}$ for each group $j =1, \cdots, p$.  The quantile estimation method doesn't perform automatic variable selection. For finding the zero vectors, i.e. the irrelevant groups of variables, hypothesis tests are required. However when model (\ref{eq1}) has a large group number $p$, it is useful to estimate simultaneously the parameter groups and to eliminate the irrelevant groups without crossing every time by a hypothesis test.  The adaptive LASSO penalties have the advantage of automatic selection and of parameter estimation (see for example \cite{Zhang-Xiang-15}, \cite{Wei-Huang-10}, \cite{Wang-Leng-08}). \\
 
 In order to introduce and study the adaptive LASSO estimator, we consider the following index set
  \[{\cal A} \equiv \{ j ; \| \eb^0_j \| \neq 0\}=\{1, \cdots, p_0\}\]
   and ${\cal A}^c \equiv \{ j ; \| \eb^0_j \| = 0\}=\{p_0+1, \cdots, p\}$ its complementary set.  The set ${\cal A}$ contains the index set corresponding to groups with nonzero true parameters.\\
  For $\eb$ a $r$-vector of  parameters, we denote by $\eb_{\cal A}$ the $r^0$-subvector of $\eb$ which contains  $\eb_j$, for $j = 1, \cdots , p_0$. Similarly, the  $(r-r^0)$-vector $\eb_{{\cal A}^c}$ contains $\eb_j$ for $j=p_0+1, \cdots, p$.\\
  In practice, the set ${\cal A}$ is unknown. Then, we must find the set ${\cal A}$ and estimate the corresponding parameters.\\
\hh In  Sections 3 and 4 we will introduce an estimator, denoted $\widehat{\eb}^*_n$, which minimizes the quantile  process penalized with an adaptive group LASSO penalty, for two cases: $p$ fixed  and   $p \rightarrow \infty$ as $n \rightarrow \infty$. We generalize the adaptive LASSO quantile estimator proposed by \cite{Ciuperca-15b} for individual variable selection to the case of group selection.  We call this estimator, adaptive group  LASSO quantile (\textit{ag\_LASSO\_Q}) estimator.\\
\hh  We say that $\widehat{\eb}^*_n$ satisfies the \textit{oracle properties} if:\\
  \textit{(i) asymptotic normality:}  $\sqrt{n}( \widehat{\eb}^*_n - \eb^0)_{\cal{A}}$ converges in law to a centred Normal distribution.\\
  \textit{(ii) sparsity property:} $\lim_{n \rightarrow \infty} \PP [{\cal A}=\{ j=1, \cdots , p ; \|\widehat{\eb}^*_{n;j} \| \neq \textbf{0}_{d_j} \} ]=1.$

  \section{Fixed $p$ case}
 In this section we propose and study the asymptotic properties of the  \textit{ag\_LASSO\_Q} estimator for the parameter $\eb$ of model (\ref{eq1}) when the group number  $p$ is fixed. \\
 We define the  \textit{ag\_LASSO\_Q} estimator by: 
 \[\widehat{\eb}^*_n \equiv \argmin_{\eb \in \R^r} Q(\eb),\]
  where $Q(\eb)$ is the penalized quantile  process with the adaptive group LASSO penalty: 
 \begin{equation}
  \label{eq3}
 Q(\eb) \equiv \sum^n_{i=1} \rho_\tau(Y_i -\sum^p_{j=1}\XX_{ij}^t \eb_j)  +\mu_n\sum^p_{j=1}  \widehat{\omega}_{n;j} \| \eb_j\|,
  \end{equation}
  with the weight  $\widehat{\omega}_{n;j} \equiv \| \widetilde{\eb}_{n;j} \| ^{- \gamma}$, $\gamma >0$.  The estimator   $\widehat{\eb}^*_n$ is written as $\widehat{\eb}^*_n =(\widehat{\eb}^*_{n;1}, \cdots, \widehat{\eb}^*_{n;p})$ and  $\widehat{\eb}^*_{n;j}$ is a  subvector of size  $d_j$, for $j=1, \cdots , p$.\\
  For a particular case of a quantile model with non-grouped variables, $d_j=1$ for all $j=1, \cdots, p$, we obtain the adaptive LASSO quantile estimator proposed and studied by \cite{Ciuperca-15b}.\\

  Before presenting the main results for $\widehat{\eb}^*_n$ in the fixed $p$ case, we give the required assumptions.\\
The tuning parameter $\mu_n$ and the constant $\gamma$ are such that, for $n \rightarrow \infty$, 
\begin{equation}
\label{cond_lambda}
\mu_n \rightarrow\infty,\quad  \frac{\mu_n}{\sqrt{n}} \rightarrow 0, \quad n^{(\gamma-1)/2} \mu_n \rightarrow \infty.
\end{equation}
 For the design $(\eeX_i)_{1 \leqslant i \leqslant n}$  we consider the following assumption:\\
\textbf{(A1)} $n^{-1} \max_{1 \leq i \leq n} \eeX_i^t \eeX_i {\underset{n \rightarrow \infty}{\longrightarrow}}  0$ and 
$n^{-1} \sum^{n}_{i=1} \eeX_i \eeX_i^t {\underset{n \rightarrow \infty}{\longrightarrow}} \eU$, with $\eU$ a $r \times r$ positive definite matrix.\\
   {For the errors $\varepsilon_i$ we suppose that:\\
 \textbf{(A2)} $(\varepsilon_i)_{1 \leqslant i \leqslant n}$ are  independent, identically distributed, with $F: {\cal B} \rightarrow [0,1]$ the distribution function and a continuous positive density $f$ in a neighborhood of $0$.   The $\tau$th quantile of $\varepsilon_i$ is zero: $\tau= F(0)$. Moreover,  
for every $e \in int({\cal B})$,  $\uu_r \in \R^r$ we have 
 \begin{equation}
 \label{rA2}
 \lim_{n \rightarrow \infty} n^{-1} \sum^n_{i=1} \int^{\xx^t_i \uu_r}_{0} \sqrt{n}[F(e+v/\sqrt{n})-F(e)] dv = \frac{1}{2} f(e) \uu^t_r 
 \eU
 \uu_r , 
 \end{equation}
 where  $\uu_r$ is the $r$-vector with all components 1. The set ${\cal B}$ is a real set, with $0 \in {\cal B}$. \\
 
 Assumption (A1) is standard for LASSO methods and (A2) is classic for quantile regression (see \cite{Ciuperca-15b},  \cite{Koenker-05}, \cite{Zou-Yuan-08}, \cite{Wu-Liu-09}). Assumption (A1) requests that the design matrix has a reasonable good behaviour.  For the tuning parameter $\mu_n$, the same conditions on (\ref{cond_lambda})  are required in  \cite{Ciuperca-15b} for adaptive LASSO quantile model but with ungrouped explanatory variables.  \\
 
 We make the remark that for ANOVA model, since in the analysis of variance there is a constraint for each level of a factor, we consider as constraint that the effect of this level is zero. Then this zero level is  not considered in the model  in order that assumption (A1) is satisfied. \\
   
In order to study the asymptotic properties of the estimator  $\widehat{\eb}^*_n$, let us consider  the index set of the groups selected by the adaptive  group  LASSO quantile method:
\[
\widehat{\cal A}^*_n \equiv \{ j \in \{1, \cdots, p\}; \| \widehat{\eb}^*_{n;j}\| \neq 0 \}  
\]
and $\widehat{A}^{*c}_n$ its complementary set.\\

The following Theorem shows that the  \textit{ag\_LASSO\_Q} estimators  with the index in the set ${\cal A}$ are asymptotically Gaussian. Then, the estimators of the nonzero parameter vectors have the same asymptotic distribution they would have if the zero parameter vectors were known. 

\begin{theorem}
\label{theorem 1}
Under assumptions (A1), (A2) and condition (\ref{cond_lambda}),  we have $\sqrt{n}( \widehat{\eb}^*_n - \eb^0)_{\cal{A}}  \overset{\cal L} {\underset{n \rightarrow \infty}{\longrightarrow}} {\cal N}\big(\textbf{0}_{r^0}, \tau(1-\tau) f^{-2}(0)\eU^{-1}_{\cal A} \big)$, with $\eU_{\cal A}$ the submatrix of $\eU$ with the row and column indices in $\{1, \cdots , d_1, d_1+1, \cdots , d_1+d_2, \cdots, \sum_{j=1}^{p_0}d_j  \}$. 
\end{theorem}

We give now the Karush-Kuhn-Tucker(KKT) optimality conditions, needed to prove the sparsity property for $\widehat{\eb}^*_n$.   \\
For all  $ j \in \widehat{\cal A}^*_n $, we have, with probability one,  the following   $d_j$ equalities  
\begin{equation}
\label{KKTi}
\tau \sum^n_{i=1} \XX_{ij} - \sum^n_{i=1} \XX_{ij} \e1_{Y_i < \eeX^t_i \widehat{\eb}^*_n}= \frac{\mu_n \widehat{\omega}_{n;j} \widehat{\eb}^*_{n;j}}{\| \widehat{\eb}_{n;j}^*\|}.
\end{equation}
For all  $ j \not \in \widehat{\cal A}^*_n $, for all $k =1, \cdots , d_j$ we have, with probability one, the following inequality
\begin{equation}
\label{KKTii}
\left| \tau \sum^n_{i=1} X_{ij,k} - \sum^n_{i=1} X_{ij,k} \e1_{Y_i < \eeX^t_i \widehat{\eb}^*_n}  \right| \leq \mu_n \widehat{\omega}_{n;j} .
\end{equation}

The following theorem shows the sparsity property of the  \textit{ag\_LASSO\_Q} estimator. This result states that the adaptive group LASSO quantile estimators of the nonzero parameter vectors are exactly nonzero with a probability converging to one when $n$ diverges to infinity. 

\begin{theorem}
\label{th_selection}
Under the assumptions of Theorem  \ref{theorem 1} and under the condition  $n^{\gamma/2-1} \lambda_n \rightarrow \infty$, as $n \rightarrow \infty$, we have $\lim_{n \rightarrow \infty }\PP [\widehat{\cal A}^*_n={\cal A}]=  1 $.
\end{theorem}

Theorem \ref{theorem 1} and Theorem \ref{th_selection} establish the asymptotic normality and the sparsity of the \textit{ag$\_$LASSO$\_$Q} estimator, which means that this estimator still share the oracle properties in the case of fixed $p$.

\begin{remark}
For the weight $\widehat{\omega}_{n;j}$ associated to the $j$th group, we considered the quantile estimator norm to the power $- \gamma$. In view of the proofs of Theorem \ref{theorem 1} and Theorem \ref{th_selection}, these two theorems remain true also when $\widetilde{\eb}_{n;j}$ is replaced by any estimator of $\eb_j$, with convergence rate $n^{-1/2}$, under assumptions (A1), (A2).
\end{remark}

\section{The case of $p$ depending on $n$}
Consider now same model (\ref{eq1}) with grouped  variables, but with the number $p$ of  groups depending on $n$: $p=p_n$ and $p_n \rightarrow \infty$ as $n \rightarrow \infty$. More precisely, we consider $p_n=O(n^c)$, with the constant $c \in (0,1)$. For readability,  we keep the notation $p$ instead of $p_n$. Similarly, we have $r=\sum^{p}_{j=1} d_j$, with $r$ depending on $n$. Always for simplicity of notation, for the design $\eeX_i$, for  the parameter $\eb$, ever if their dimension depends on $n$, we do not put subscript $n$. \\

We will first find the  convergence rate of the  quantile estimator  $\widetilde{\eb_n} $ of $\eb$. Afterwards, we will  propose for $\eb$ an adaptive group  LASSO quantile estimator. Even though the number $p$ diverges as $n \rightarrow \infty$,  this estimator keeps the   oracle properties. \\

Since the design size depends on $n$, we need reconsider  the assumptions on $\eeX_i$. Then, let us consider the following assumptions for the errors $(\varepsilon_i)$,  design $(\eeX_i )$ and for the number  $p$ of groups:\\
\textbf{(A3)}  $(\varepsilon_i)_{1 \leq i \leq n}$ are i.i.d. Let  $F$ be   the distribution function  and  $f$ be the density function of $(\varepsilon_i)$. The density function $f$ is continuously,  strictly positive in a neighbourhood of zero and has a bounded first derivative in the neighbourhood of 0. The $\tau$th quantile of $\varepsilon_i$ is zero: $\tau= F(0)$. \\
\textbf{(A4)} There exist two constants $0 < m_0 \leq M_0 < \infty$, such that  $ m_0 \leq \lambda_{\min} ( {n}^{-1} \sum^n_{i=1} \eeX_i \eeX_i^t) \leq \lambda_{\max} ( {n}^{-1} \sum^n_{i=1} \eeX_i \eeX_i^t) \leq M_0 $. \\
\textbf{(A5)}  $\pth{ {p}/{n}}^{1/2} \max_{1 \leqslant i \leqslant n} \| \eeX_i\| \rightarrow 0$, as $n \rightarrow \infty$.\\
\textbf{(A6)} $p$ is such that $p=O(n^c)$, with $0< c < 1$.\\

Since   $p \rightarrow \infty$, condition (\ref{rA2}) of assumption (A2) for the case $p$ fixed is now  replaced by  $f'$ bounded in the neighborhood of 0. This assumption also  been considered for always high-dimensional quantile model, with seamless $L_0$ penalty by \cite{Ciuperca-15}. In \cite{Ciuperca-15}, assumptions (A4) and (A5) are also required.  Assumption (A6) was considered by \cite{Zhang-Xiang-15} for an high-dimensional linear model where the objective function is   the error squares sum, penalized with an adaptive group LASSO penalty. Assumptions (A4), (A5), (A6) are also required for an high-dimensional linear model by \cite{Zou-Zhang-09}, which penalize the LS process with adaptive elastic-net penalty.  In respect to the case $p$ fixed, assumptions (A4) and (A5) are the similar of  (A1). \\

We will start by finding the convergence rate of  quantile estimator (\ref{eq2}) in the case $p \rightarrow \infty$ as $n \rightarrow \infty$.  For this, consider the  quantile process:
\[
G_n(\eb) \equiv \sum^n_{i=1} \rho_\tau(Y_i-\eeX^t_i \eb).
\]
For the quantile estimator existence, we assume that the total number $r$  of parameters is strictly less than $n$. \\

We recall that in the case  $p$ fixed, the  convergence rate of the  quantile estimator $ \widetilde{\eb}_n$ is of order  $n^{-1/2}$ (see for example \cite{Koenker-05}).
We will show that, the quantile estimator has the convergence of order $(p/n)^{1/2}$, when the explanatory group variable number diverges with the sample size. In view of the proof of Lemma \ref{Lemma 2.1}, the convergence rate of $ \widetilde{\eb}_n$ depends only of $p$ and not of total number $r$ of parameters, thanks to assumption (A5).  One needs the convergence rate of the quantile estimator is necessary for  studying the asymptotic behaviour of the penalty  which intervenes  in adaptive group LASSO quantile process. 

\begin{lemma}
\label{Lemma 2.1}
Under assumptions (A3)-(A6), we have $\| \widetilde{\eb}_n- \ebo\|=O_{\PP}\pth{\sqrt{\frac{p}{n}}}$.
\end{lemma}

Consider now the following adaptive group  LASSO  quantile (\textit{ag\_LASSO\_Q}) estimator:
\[\widehat{\eb}^*_n \equiv \argmin_{\eb \in \R^d} \pth{ \frac{1}{n} G_n(\eb)+\lambda_n \sum^p_{j=1} \widehat{\omega}_{n;j} \| \eb_j\| },
\]
where $\lambda_n$ is a tuning parameter  (positive) and the weights of the LASSO penalty are   $\widehat{\omega}_{n;j} \equiv \| \widetilde{\eb}_{n;j} \| ^{- \gamma}$, with $\gamma >0$. The relation between the tuning parameter   $\mu_n$ of relation (\ref{eq3}) for the case  $p$ fixed  and  $\lambda_n$ for the case $p$ depending on $n$ is $\lambda_n = \mu_n /n$. We prefer to consider these forms  as tuning parameter and as objective  process,   for having a similarity with the  adaptive group LASSO  LS (\textit{ag\_LASSO\_LS})  case considered by  \cite{Zhang-Xiang-15}. \\

In order to study the asymptotic normality of $ \widehat{\eb}^*_n$ we need to impose an additional condition on the total number  of nonzero parameters. More precisely, $r^0$ it is assumed to be the same order as $p_0$. This is for controlling the penalty,  so that it is smaller than the quantile process.\\

Concerning the size of the nonzero  parameter vectors, we take the following assumption: \\
\textbf{(A7)} $r^0=O(p_0)$.\\
For the smallest nonzero vector norm and on constant $c$ of assumption (A6) we assume:\\
\textbf{(A8)} Let us denote $h_0 \equiv \min_{1 \leqslant j \leqslant p_0} \| \eb^0_j\|$. There  exists a  constant $M>0$ such that $M \leq n^{- \alpha} h_0$ and $\alpha > (c-1)/{2}$. \\

These two assumptions  were also found in the paper \cite{Zhang-Xiang-15}, for  \textit{ag\_LASSO\_LS} method in high-dimensional linear model, but with a supplementary condition for $r$: $r=O(p)$. Here, we do not need this requirement, since assumption (A5) is imposed. On the other hand, in \cite{Zhang-Xiang-15}, instead of assumption (A5) the condition $n^{-1/2} \max_{1 \leqslant i \leqslant n}  \| \eeX_{i,{\cal A}} \|^2 \rightarrow 0 $, as $n \rightarrow \infty$, is required. \\

The following theorem gives the convergence rate of the  \textit{ag\_LASSO\_Q} estimator when  $p \rightarrow \infty$. We obtain the same convergence rate  that of  quantile estimator when group number diverges. This convergence rate is also obtained by \cite{Zhang-Xiang-15} for the \textit{ag\_LASSO\_LS} estimator, but for errors $(\varepsilon)_{1 \leqslant i \leqslant n}$ with  mean zero and   bounded variance.

\begin{theorem}
\label{th_vconv}
Under assumptions (A3)-(A6), (A8) and the  tuning parameter $(\lambda_n)_{n \in \N}$ satisfying $\lambda_n n^{(1+c)/{2} - \alpha \gamma} \rightarrow 0$, as $n \rightarrow \infty$, we have   $\| \widehat{\eb}^*_n- \ebo\|=O_{\PP}\pth{\sqrt{\frac{p}{n}}}$.
\end{theorem}

The following theorem shows the  oracle properties for \textit{ag\_LASSO\_Q} estimator when the number $p$ of groups diverges. We denote by  $\eeX_{i,{\cal A}}$ a $r^0$-vector which contains the sub-vectors $\XX_{i,j}$, for $j \in \{1, \cdots , p_0 \}$. 
\begin{theorem}
\label{Theorem 2SPL} Suppose that assumptions (A3)-(A6), (A8) are  satisfied  and also that the tuning parameter satisfies   $\lambda_n n^{(1-c)  {(1+\gamma)}/{2}} \rightarrow \infty$, $\lambda_n n^{(c+1)/{2} - \alpha \gamma} \rightarrow 0$, as $n \rightarrow \infty$.  Then:\\
(i) $\PP \cro{\widehat{\cal A}^*_n={\cal A}}\rightarrow 1$,  for $n \rightarrow \infty $.\\
(ii) If moreover assumption (A7) holds, then, for any  vector $\eu$ of size $r^0$ such that $\| \eu\|=1$, with notation $\eU_{n,{\cal A}} \equiv n^{-1} \sum^n_{i=1} \eeX_{i,{\cal A}}  \eeX_{i,{\cal A}}^t$, we have  $\sqrt{n} (\eu^t \eU^{-1}_{n,{\cal A}} \eu)^{-1/2} \eu^t ( \widehat{\eb}^*_n - \eb^0)_{\cal{A}}  \overset{\cal L} {\underset{n \rightarrow \infty}{\longrightarrow}} {\cal N}\big(0, \tau (1- \tau ) f^{-2}(0) \big)$.
\end{theorem}

For the  tuning parameter $\lambda_n$, the same conditions are required in \cite{Zhang-Xiang-15} such that, the \textit{ag\_LASSO\_LS} estimator  in an high-dimensional linear model satisfies the oracle properties. 

\begin{remark}
As for the case  $p$ fixed, we considered the weight  $\widehat{\omega}_{n;j} = \| \widetilde{\eb}_{n;j} \| ^{- \gamma}$, with $\widetilde{\eb}_{n;j}$ the quantile estimator of the $d_j$-vector $\eb_j$, for any $j=1, \cdots, p$. In view of the proof of Theorem \ref{Theorem 2SPL}, the oracle properties for \textit{ag\_LASSO\_Q} estimator remain true also when $ \widetilde{\eb}_{n;j}$ is replaced by any $(p/n)^{1/2}$-estimator of $\eb_j$, under assumptions (A3)-(A6).
\end{remark}

\begin{remark}
If $h_0$, defined in assumption (A8), doesn't depend on $n$, then $\alpha=0$. In this case, the conditions required on $(\lambda_n)_{n\in \N}$ in Theorem \ref{th_vconv} imply $\gamma >  {2 c }/{(1-c)}$, and then $\gamma$ can take values in the interval $ (0, \infty)$. The value of $\gamma$ increase with that of $c \in (0,1)$. For example, if $c={1}/{2}$ then $\gamma >2$. 
\end{remark}

\section{Simulations}
In order to evaluate the performance of the proposed estimation method, Monte Carlo simulations are realized in this section. To assess this performance we compare  the \textit{ag\_LASSO\_Q} and \textit{ag\_LASSO\_LS} estimation methods.\\

The design $\eeX_i$ is generated in the same way as in  paper \cite{Wei-Huang-10}: $\eeX=(\XX_1, \cdots , \XX_p)$, with the group explanatory variables  $\XX_j=(X_{5(j-1)+1}, \cdots, X_{5j})$, for all $j=1, \cdots, p$. We first generate $r=5p$  independent random  variables   $R_1, \cdots , R_r$ of standard normal distribution. We also generate the  variables $Z_j$  of multivariate normal distribution with mean zero and covariance $Cov(Z_{j_1}, Z_{j_2})=0.9^{| j_1-j_2|}$. Finally, the  variables $X_1, \cdots , X_{r}$ are generated as:
\[
X_{5(j-1)+k}=\frac{Z_j+R_{5(j-1)+k}}{\sqrt{2}}, \qquad 1 \leq j \leq p, \quad 1 \leq k \leq 5.
\]
Two model errors are considered: Normal ${\cal N}(0,3^2)$ and Cauchy ${\cal C}(0,3^2)$. For the parameters we take: $\eb^0_1=(0.5,1,1.5,1,0.5)$, $\eb^0_2=(1,1,1,1,1)$, $\eb^0_3=(-1,0,1,2,1.5)$, $\eb^0_4=(-1.5,1,0.5,0.5,0.5)$ and all other parameters are zero vectors. The nonzero vectors were also considered in Example 2 of \cite{Wei-Huang-10} for errors ${\cal N}(0,3^2)$, $p=10$, when the parameters were  estimated by LS method with  adaptive group   LASSO penalty. \\
The constant $c$ of assumption (A6) is $c=0.43$. Then,  we will consider the following value couples for $n$ and $p$: $(30,5)$, $(60,5)$, $(60,10)$, $(100,10)$, $(200,10)$, $(400,15)$, $(1000,25)$ and $(1000,100)$. On the other hand, $p^0$  will always be equal to 4. The response  variable   $Y$ is generated as: $Y_i=\sum^p_{j=1}\XX_{ij}^t \eb^0_j+\varepsilon_i$, for $i=1, \cdots, n$.\\
We will compare the obtained results by the adaptive group  LASSO quantile method,  proposed in this paper,  with those obtained by  the adaptive group  LASSO LS method, proposed by  \cite{Wei-Huang-10}, \cite{Zhang-Xiang-15}. \\
For simulations, we used the R language. After a scale transformation, we can use the group LASSO methods instead of the adaptive  LASSO group methods. Then, in order to calculate the adaptive group LASSO  LS estimations we have used the function \textit{grpreg} of package \textit{grpreg}, the tuning parameter being chosen on a value  grid, using the  AIC criterion. In order to  calculate the adaptive group  LASSO quantile estimations, we have used the  function \textit{groupQICD} of package \textit{rqPen} and the  tuning parameter varies on a value  grid.\\
For each considered case,  1000  Monte Carlo replications was made.\\

In Table  \ref{Tabl1} we give how the two estimation methods  identify the parameter vectors (zero or nonzero), for the part that contains the  four nonzero parameter vectors $\eb^0_j$, $j=1, \cdots, 4$,   and for the part with  $p-4$ zero vectors. We present the minimum, three quartiles (Q1, median, Q3), the mean and the  maximum of the  number of nonzero vectors ($j=1, \cdots, 4$), respectively, zero ($j=5, \cdots, p$), found by the two estimation methods.\\
For $n$ large (equal to 100, 200, 400, 1000), we observe that for errors of ${\cal N}(0,3^2)$  law, the two estimation methods well  identify  the zero and nonzero   parameter vectors.  However, for Cauchy errors, the  \textit{ag\_LASSO\_LS} method  poorly identifies nonzero vectors (the  group of the four significant variables). The zero vectors are  very well identified by the two methods.\\
For $n$ small (equal to 30 or 60), the two estimation  methods well identify the four relevant variable groups, that errors are Normal or Cauchy (except for  \textit{ag\_LASSO\_Q}, in the case $n=60$, $p=5$, $\varepsilon \sim {\cal C}(0,3^2)$). However, the $(p-4)$  irrelevant variable groups are not well identified by the \textit{ag\_LASSO\_LS} method. \\

\textit{Conclusion}\\
\hh For gaussian errors, the  \textit{ag\_LASSO\_LS} method  identifies well the two (relevant and irrelevant) variable groups for  $n$ large. For  $n$ small, the irrelevant variable groups are not well identified.  For  Cauchy errors, this method, either does not identify the relevant variable groups or irrelevant variable groups, regardless of the value $n$. Then, for Cauchy errors, the \textit{ag\_LASSO\_LS} estimations do not have the sparsity property. \\
\hh The \textit{ag\_LASSO\_Q} method, for the two  types of  errors, identifies the two    variable groups (significant and irrelevant), the precision increasing with $n$. Then, the \textit{ag\_LASSO\_Q} estimations have the sparsity property. \\
\hh We conclude then that the simulations confirm the theoretical results for the \textit{ag\_LASSO\_Q} estimators.
\begin{table}[p] 
{\scriptsize
\caption{\footnotesize Model selection results by \textit{ag\_LASSO\_Q} and \textit{ag\_LASSO\_LS} methods for  $p^0=4$, $c=0.43$, errors   ${\cal N}(0,3^2)$ or ${\cal C} (0,3^2)$. }
\begin{center}
\rotatebox{90}{
\begin{tabular}{|c|c|c|c|c|c|c|c|c|c|c|c|c|c|c||c|c|c|c|c|c|c|c|c|c|c|c|c|}\hline
 & & &  \multicolumn{12}{c||}{nonzero true parameter vectors} & \multicolumn{13}{c|}{zero true parameter vectors} \\
  \cline{4-28} 
 n &p & $\varepsilon$ & \multicolumn{2}{c|}{min} & \multicolumn{2}{c|}{Q1} & \multicolumn{2}{c|}{\textbf{Q2}} & \multicolumn{2}{c|}{mean} & \multicolumn{2}{c|}{Q3} & \multicolumn{2}{c||}{max} & \textbf{true} & \multicolumn{2}{c|}{min} & \multicolumn{2}{c|}{Q1} & \multicolumn{2}{c|}{\textbf{Q2}} & \multicolumn{2}{c|}{mean} & \multicolumn{2}{c|}{Q3} & \multicolumn{2}{c|}{max} \\
   & & & LS & Q & LS & Q & LS & Q & LS & Q & LS & Q & LS & Q & & LS & Q & LS & Q & LS & Q & LS & Q & LS & Q & LS & Q \\ \hline
    30      & 5 & {\cal C} &0 & 0 & 4 & 3 & \textbf{4}&  \textbf{3} & 3.5 & 3.3 &4 & 4 &  4 & 4 &\textbf{1} & 0 & 0 & 0 & 0 & \textbf{0} & \textbf{1} &0.1 &0.7 & 0 &  1 & 1 & 1\\
     & & {\cal N} & 4& 3 & 4 & 4 & \textbf{4}&\textbf{4} & 4 & 3.8 & 4 &4 & 4 &4 &\textbf{1} & 0 & 0 & 0 & 0 & \textbf{0} &\textbf{1}  & 0.01 & 0.5 &0 & 1& 1 & 1  \\ \hline 
    60    & 5 & {\cal C} &0 & 1 & 0 & 3 &\textbf{1} &\textbf{4} & 1.3 & 3.6 &  4& 4& 4 & 4 &\textbf{1}& 0 & 0 & 0 & 1 & \textbf{1} & \textbf{1} & 0.7 & 0.8 & 1& 1& 1 & 1 \\
     & & {\cal N} &4 & 3 & 4 & 4 & \textbf{4}& \textbf{4}& 4 & 3.9 & 4 & 4& 4 &4 &\textbf{1} & 0 & 0 & 0 &  0& \textbf{ 0}& \textbf{1} & 0.1 & 0.7 &0 &1 & 1 & 1  \\ \hline 
    60       & 10 & {\cal C} & 0&  0& 4 &3 & \textbf{4}&\textbf{4} & 3.4 &3.5  &4  & 4& 4 &4 &\textbf{6} & 0 & 2 & 0 & 5 & \textbf{0} & \textbf{5} & 0.8 & 5.2 & 0 & 6 & 6 & 6 \\
   & & {\cal N} &4 & 3 & 4 & 4 & \textbf{4}&\textbf{4} & 4 & 3.9 &4  & 4& 4 & 4 &\textbf{6}& 0 &1  & 0 & 4 &\textbf{0} & \textbf{4} & 0.06 & 4.3 & 0&5 & 6 &  6 \\ \hline 
     100    & 10 & {\cal C} &0 & 3 & 0 & 4 & \textbf{0}&\textbf{4} & 0.9 & 3.9 & 1 &4&4  & 4 &\textbf{6}& 0 & 2 & 5 & 4 &\textbf{6}  & \textbf{5} &5  &  5 & 6& 6 & 6 &6\\
      & & {\cal N} &4 & 3 & 4 & 4 & \textbf{4}&\textbf{4} & 4 & 3.9 & 4 &4 &4&4 &\textbf{6}& 0 &1 & 0 & 4 & \textbf{0} &\textbf{ 5} & 2 & 4.4 &  5& 5 &6 &  6  \\ \hline 
   200 & 10 & {\cal C} & 0& 3 & 0 & 3 & \textbf{0}& \textbf{4} & 0.5 & 3.7 & 1 & 4& 4 & 4 &\textbf{6}& 0 & 2 & 6 & 6 & \textbf{6} & \textbf{6} & 5.7 & 5.7 & 6& 6&  6& 6 \\
      & & {\cal N} & 4& 3 & 4 & 4 &\textbf{ 4}& \textbf{4} & 4  & 3.9 & 4 & 4& 4 & 4 &\textbf{6}& 0 & 3 & 0 & 5 & \textbf{4} & \textbf{6} & 3.1 & 5.5 & 6& 6&  6& 6 \\ \hline 
      400 & 15 & {\cal C} & 0& 3 & 0 & 3 & \textbf{0}& \textbf{4} & 0.4 & 3.6 &\textbf{11} & 0 & 4& 4 & 4& 8 & 10 & 11 & 11 & \textbf{11} & \textbf{11} & 10.8 & 10.9 & 11& 11&  11& 11 \\
      & & {\cal N} & 4& 3 & 4 & 4 & \textbf{4}& \textbf{4} & 4  & 3.99 & 4 & 4& 4 & 4 &\textbf{11}& 0 & 9 & 9 & 11 & \textbf{10} & \textbf{11} & 9.3 & 10.8 & 11& 11&  11& 11 \\ \hline 
   1000 & 25 & {\cal C} & 0& 4 & 0 & 4 & \textbf{0}& \textbf{4} & 0.5 & 4 & 1 & 4& 4 & 4 &\textbf{21}& 19 & 19 & 21 & 21 & \textbf{21} & \textbf{21} & 20.8 & 20.7 & 21& 21&  21 & 21 \\
     & & {\cal N} & 4& 4 & 4 & 4 & \textbf{4}& \textbf{4} & 4  & 4 & 4 & 4& 4 & 4 &\textbf{21}& 17 & 20 & 20 & 21& \textbf{21} & \textbf{21} & 20.3 & 20.9 & 21& 21&  21& 21  \\ 
   \cline{2-27} 
   & 100 & {\cal C} & 0& 3 & 0 & 3 & \textbf{ 0} & \textbf{4} & 0.3 & 3.5 & 0& 4& 3&  4 &\textbf{96}&94 & 96 & 96 & 96 & \textbf{96} & \textbf{96} & 95.8 & 96 & 96 & 96 &96 &  96 \\
      & & {\cal N} & 4& 3 & 4 & 4 & \textbf{4}& \textbf{4} & 4  & 3.98 & 4 & 4& 4 & 4 &\textbf{96} & 92 & 95 & 95 & 96 & \textbf{96} & \textbf{96} & 95.4 & 96 & 96& 96&  96& 96  \\ \hline
      \end{tabular}
}
 \end{center}
\label{Tabl1} 
}
\end{table}

\section{Proofs}
In this section we provide the proofs of  all results presented in Sections 3 and 4.
\subsection{Proofs for results of Section 3}
\noindent {\bf Proof of Theorem \ref{theorem 1}}. The proof is similar to that of  Theorem 4.1 of \cite{Zou-Yuan-08}.\\
We denote $\sqrt{n}( \widehat{\eb}^*_n - \eb^0) \equiv \widehat{\eu}_n$, and in general $\sqrt{n}(\eb- \ebo)\equiv \eu\equiv (\eu_1, \cdots, \eu_p) $, with $\eu_j=(u_{j,1}, \cdots , u_{j,d_j})$, for $j=1, \cdots , p$.\\
Since $Y_i=\eeX_i^t \ebo+ \varepsilon_i$, then 
$
Y_i-\eeX_i^t \eb=\frac{\eeX_i^t \eu}{\sqrt{n}} +\varepsilon_i.$ 
Let us consider the following random variables
\begin{eqnarray}
\label{Di}
{\cal D}_i & \equiv & (1-\tau) \e1_{\varepsilon_i <0}- \tau \e1_{\varepsilon_i \geq 0},  \\
v_n & \equiv &  \frac{1}{\sqrt{n}} \sum^n_{i=1} {\cal D}_i, \nonumber \\ 
B_n(\eu) & \equiv & \sum^n_{i=1} \int^{\eeX^t_i \eu/\sqrt{n}}_0 [\e1_{\varepsilon_i < t} - \e1_{\varepsilon_i < 0} ]dt
 \nonumber 
\end{eqnarray}
and the random vector  
\[
\bfz_n  \equiv     \frac{1}{\sqrt{n}} \sum^n_{i=1} \eeX_i  {\cal D}_i. 
  \]
Obviously, $\eE[{\cal D}_i]=0$ et $\eE[\bfz_n]=\textbf{0}_r$.
By the CLT, using assumptions (A1) and (A2), we have
\begin{equation}
\label{eq4}
 \bfz_n  \overset{\cal L} {\underset{n \rightarrow \infty}{\longrightarrow}} {\cal N}(\textbf{0}_r, \tau(1-\tau)\eU ) , \qquad 
v_n  \overset{\cal L} {\underset{n \rightarrow \infty}{\longrightarrow}} {\cal N}(0, \tau(1-\tau)) .
\end{equation}
The vector $\widehat{\eu}_n$ is the minimizer of the following random process:
\[
L_n(\eu) \equiv \sum^n_{i=1} \big[ \rho_\tau \big( \varepsilon_i - \frac{\eeX^t_i \eu}{\sqrt{n}}\big)- \rho_\tau(\varepsilon_i) \big]+\mu_n\sum^p_{j=1}  \widehat{\omega}_{n;j} \big[ \left\|\eb^0_j +\frac{\eu_j}{\sqrt{n}} \right\| - \| \eb^0_j \|  \big],
\]
which  can be written under the following form:
\begin{equation}
\label{Lnu}
L_n(\eu)=[\bfz^t_n \eu+B_n(\eu)]+ \mu_n\sum^p_{j=1}  \widehat{\omega}_{n;j}  \big[ \left\|\eb^0_j +\frac{\eu_j}{\sqrt{n}} \right\| - \| \eb^0_j \|  \big] \frac{\sqrt{n}}{\sqrt{n}}.
\end{equation}
We first study  the last sum of the right hand side of (\ref{Lnu}). \\
\underline{For all  $j \leq p_0$} (thus $\|\eb^0_j\| \neq 0$) we have, since the quantile estimators are consistent:
\begin{equation}
\label{sl1}
\widehat{\omega}_{n;j} \overset{\PP} {\underset{n \rightarrow \infty}{\longrightarrow}}  \| \eb^0_j \|^{- \gamma} \neq 0
 \end{equation}
and by elementary calculus
 \begin{equation}
\label{sl2}
 \sqrt{n} \big[\| \eb^0_j+n^{-1/2} \eu_j \| - \|\eb^0_j\|  \big]  {\underset{n \rightarrow \infty}{\longrightarrow}}  \frac{\eu^t_j \eb^0_j}{\| \eb^0_j\|}.
 \end{equation}
 Then, using condition $\mu_n n^{-1/2} \rightarrow 0$, when $n \rightarrow \infty$, taking into account relations  (\ref{sl1}}) and (\ref{sl2}), by  Slutsky's Lemma, we have:
 \begin{equation}
 \label{cas1}
 \mu_n \sum^{p_0}_{j=1} \widehat{\omega}_{n;j}  \big[\| \eb^0_j+n^{-1/2} \eu_j \| - \|\eb^0_j\|  \big] \overset{\PP} {\underset{n \rightarrow  \infty}{\longrightarrow}} 0.
 \end{equation}
  \underline{For $j >p_0$}, we have $\eb^0_j =\textbf{0}_{d_j}$. Then: $\sqrt{n} \big[\| \eb^0_j+n^{-1/2} \eu_j \| - \|\eb^0_j\|  \big] =\|\eu_j\|$. Since $\widehat{\omega}_{n;j}  \overset{\PP} {\underset{n \rightarrow \infty}{\longrightarrow}}  \infty$,  by assumption $n^{(\gamma -1)/2} \mu_n \rightarrow \infty$, we have $n^{-1/2} \mu_n \widehat{\omega}_{n;j}  \overset{\PP} {\underset{n \rightarrow \infty}{\longrightarrow}}  \infty $. Thus
 \begin{equation}
 \label{cas2}
 \mu_n \widehat{\omega}_{n,j} \frac{\sqrt{n}}{\sqrt{n}} \big[\| \eb^0_j+n^{-1/2} \eu_j \| - \|\eb^0_j\|  \big]  \overset{\PP} {\underset{n \rightarrow \infty}{\longrightarrow}}  \left\{
 \begin{array}{lll}
 0, & \textrm {if } & \eu_j=\textbf{0}_{d_j} \\
  \infty, & \textrm {if } & \|\eu_j\| \neq 0.
 \end{array}
 \right.
 \end{equation}
Then, taking into account relations  (\ref{cas1}) and (\ref{cas2}), we have the following result for the  third term of the right hand side of (\ref{Lnu}):  
\begin{equation}
\label{eq1112}
\mu_n\sum^p_{j=1}  \widehat{\omega}_{n;j}  \big[ \left\|\eb^0_j +\frac{\eu_j}{\sqrt{n}} \right\| - \| \eb^0_j \|  \big] \frac{\sqrt{n}}{\sqrt{n}}  \overset{\PP} {\underset{n \rightarrow \infty}{\longrightarrow}} \sum^p_{j=1}W(\eb^0_j,\eu),
\end{equation}
with
 \[
 W(\eb^0_j,\eu)\equiv \left\{
 \begin{array}{lll}
 0, &\textrm {if } & \eb^0_j \neq \textbf{0}_{d_j} \\
 0,& \textrm { if } & \eb^0_j = \textbf{0}_{d_j}  \textrm{ and } \eu_j =\textbf{0}_{d_j}\\
 \infty , & \textrm { if } & \eb^0_j = \textbf{0}_{d_j}  \textrm{ and } \eu_j \neq \textbf{0}_{d_j}.
 \end{array}
 \right.
 \]
On the other hand, by the two results of (\ref{eq4}), we have for the first two terms  of the right hand side of (\ref{Lnu}), with $\bfz $ a random $d$-vector of law ${\cal N}(\textbf{0}_d, \tau(1-\tau)\eU )$, that
 \[
   \bfz_n^t \eu \overset{\cal L} {\underset{n \rightarrow \infty}{\longrightarrow}} \bfz^t \eu, \qquad 
 B_n(\eu) \overset{\PP} {\underset{n \rightarrow  \infty}{\longrightarrow}}  \frac{1}{2} f(0) \eu^t \eU \eu.
 \]
Taking into account  these last  two results and relation  (\ref{eq1112}), then, $L_n(\eu)$ of relation  (\ref{Lnu}) has an asymptotic distribution:
 \[
 L_n(\eu) \overset{\cal L} {\underset{n \rightarrow \infty}{\longrightarrow}} \bfz^t \eu +  \frac{1}{2} f(0) \eu^t \eU \eu +\sum^p_{j=1}W(\eb^0_j,\eu).
 \]
Let us denote  $\eu=(\eu_1,\eu_2)$ with $\eu_1$ of size $r^0$, $\eu_2$ of size  $r-r^0$ and $\widehat{\eu}_n =(\widehat{\eu}_{1n}, \widehat{\eu}_{2n})$, where $\widehat{\eu}_{1n}$ contains the first $\sum_{j=1}^{p_0}d_j=r^0$ elements of $\eu$. Since  $\widehat{\eu}_n =\argmin_{\eu \in \R^r} L_n(\eu)$, we obtain that  $\widehat{\eu}_{2n} \overset{\PP} {\underset{n \rightarrow \infty}{\longrightarrow}} \textbf{0}_{r-r^0}$ et $ \widehat{\eu}_{1n} \overset{\cal L} {\underset{n \rightarrow \infty}{\longrightarrow}} {\cal N}(\textbf{0}_{r^0}, \tau(1-\tau) f^{-2}(0)\eU^{-1}_{\cal A})$. 
\hspace*{\fill}$\blacksquare$ \\

\noindent {\bf Proof of Theorem \ref{th_selection}}.
By  Theorem \ref{theorem 1},  for all $ j \in {\cal A}$  we have that  $\sqrt{n} (\widehat{\eb}^*_{n;j} - \eb^0_j)\overset{\cal L} {\longrightarrow} {\cal N}(\textbf{0}_{d_j}, \tau(1- \tau) f^{-2}(0) \eU_{{\cal A}_j} ) $ as $n \rightarrow \infty$. The square matrix   $\eU_{{\cal A}_j}$   of size $d_j \times d_j$ is  the submatrix of $\eU$ with the row and column indices in $\{ d_{j-1}+1, d_{j-1}+2, \cdots, d_j \}$, with $d_0=0$. Since $\eb^0_j \neq \textbf{0}_{d_j}$, then
\begin{equation}
\label{ecA1}
\lim_{n \rightarrow \infty}\PP [{\cal A} \subseteq \widehat{{\cal A}}^*_n] = 1 .
\end{equation}
To finish the proof we show that for all $j \not \in {\cal A}$ we have $\PP[j \in \widehat{\cal A}^*_n] \rightarrow 0$ as $n \rightarrow \infty$. Since $j \not \in {\cal A} $, then $\eb^0_j=\textbf{0}_{d_j}$. 
Considering the Eucildean norm for  equalities (\ref{KKTi}) we have with probability one, since we suppose $j \in \widehat{\cal A}^*_n$, that:
\[
\mu_n \widehat{\omega}_{n;j} < 2 \| \sum^n_{i=1} \XX_{ij} \| \leq 2 \sum^n_{i=1}  \|  \XX_{ij}\| = 2 \sum^n_{i=1} \pth{ \sum^{d_j}_{k=1}X_{ij,k}^2 }^{1/2}.
\]
By  the  Cauchy-Schwarz inequality, we have that,
\[
\sum^n_{i=1}\frac{1}{n} \pth{ \sum^{d_j}_{k=1}X_{ij,k}^2 }^{1/2} \leq \pth{\frac{1}{n} \sum^n_{i=1} \pth{\sum^{d_j}_{k=1}X_{ij,k}^2}}^{1/2} = \pth{\frac{1}{n} \sum^n_{i=1} \| \XX_{ij} \|^2}^{1/2}.
\]
Then, taking into account assumption (A1),  there exists a bounded constant $C_1 >0$ such that
\begin{equation}
\label{eq5}
\frac{1}{n} \mu_n \widehat{\omega}_{n;j} < 2 \pth{\frac{1}{n} \sum^n_{i=1} \| \XX_{ij} \|^2}^{1/2} \leq C_1 < \infty.
\end{equation}
On the other hand,   left-hand side of inequality (\ref{eq5}), can be written:
\[
\frac{{\mu_n \widehat{\omega}_{n;j}}}{n}=\frac{\mu_n}{n^{\gamma/2} \|\widetilde{\eb}_{n;j} \|^\gamma}  \frac{n^{\gamma /2}}{n}.
\]
Since we have supposed that  $j \in \widehat{\cal A}^*_n$ and $j \not \in {\cal A}$, we have that for all  $\epsilon>0$, there exists $\eta_\epsilon >0$ such that $\PP[n^{-1/2} \|\widetilde{\eb}_{n;j}\|^{-1} >\eta_\epsilon]>1-\epsilon$. The last two relations, together with the supposition $n^{\gamma/2-1} \mu_n \rightarrow \infty$,  imply, for all constant $A>0$,
\begin{equation}
\label{eq6}
\lim_{n \rightarrow \infty}\PP\bigg[\frac{{\mu_n \widehat{\omega}_{n;j}}}{n}>A \bigg]=1 .
\end{equation}
Then,    relations (\ref{eq5}) and (\ref{eq6}) are in contradiction. Thus
\begin{equation}
\label{ecA2}
\lim_{n \rightarrow \infty}\PP[j \in {\cal A}^c \cap \widehat{\cal A}^*_n] = 0.
\end{equation}
The theorem follows from relations (\ref{ecA1}) and (\ref{ecA2}). 
\hspace*{\fill}$\blacksquare$ \\

\subsection{Proofs for results of Section 4}

\noindent {\bf Proof of Lemma \ref{Lemma 2.1}}.
We show that for all $ \epsilon > 0$, there exists a  constant  $B_\epsilon >0$ (without loss of generality, we consider $B_\epsilon >0$, otherwise we take $|B_\epsilon|$) large enough such that for $n$ large enough:
\begin{equation}
\label{eq8}
\PP \cro{ \inf_{\| \eu \| =1} G_n\pth{\ebo+B_\epsilon \sqrt{\frac{p}{n}} \eu}>G_n(\ebo)} \geq 1-\epsilon .
\end{equation}
For this, we consider for some constant $C>0$, the expectation of the difference:
\[
\eE \cro{ G_n\pth{\ebo+C \sqrt{\frac{p}{n}} \eu} - G_n(\ebo) } = \sum^n_{i=1} \eE \cro{\rho_\tau \pth{\varepsilon_i - C \sqrt{\frac{p}{n}} \eeX^t_i \eu} - \rho_\tau(\varepsilon_i) }
\]
\begin{equation}
\label{ee}
 \qquad \qquad  =\sum^n_{i=1} \eE \cro{\int^{C \sqrt{\frac{p}{n}} \eeX^t_i \eu} _0  \e1_{0 < \varepsilon_i < t } dt}=\sum^n_{i=1} \int^{C \sqrt{\frac{p}{n}} \eeX^t_i \eu} _0 \cro{F (t)- F (0)}dt.
\end{equation}
By assumption (A6) we have  $p/n \rightarrow 0$. Moreover, by assumption (A5), 
we have that $\sqrt{\frac{p}{n}} \eeX_i^t \eu=o(1)$, for $\| \eu \|=1$. Thus, by mean value theorem and since the density $f$ has a bounded first derivative in the neighbourhood of 0,  relation (\ref{ee}) becomes:
\[
    =\frac{f(0)}{2} C^2\frac{p}{n} \sum^n_{i=1} (\eeX^t_i \eu)^2+ o\pth{\frac{p}{n} \sum^n_{i=1} \eu^t (\eeX^t_i \eeX_i)\eu}. \qquad \qquad
\]
Then, taking into account assumption (A4),
\begin{equation}
\label{eq9}
\frac{1}{n}\eE \cro{ G_n\pth{\ebo+C \sqrt{\frac{p}{n}} \eu} - G_n(\ebo) } = C^2 \frac{f(0)}{2} \frac{p}{n} \frac{1}{n} \sum^n_{i=1} (\eeX^t_i \eu)^2 (1+o(1)).
\end{equation}
Let be the following random  variable 
\[
{\cal R}_i  \equiv  \rho_\tau\left(\varepsilon_i -C \sqrt{\frac{p}{n}} \eeX^t_i \eu\right) -\rho_\tau(\varepsilon_i)-C\sqrt{\frac{p}{n}}  {\cal D}_i \eeX^t_i \eu \nonumber 
\]
and the following random  vector  
\[
\textbf{W}_n \equiv C\sqrt{\frac{p}{n}} \sum^n_{i=1}  {\cal D}_i \eeX^t_i,
\]
with the random  variable    ${\cal D}_i$ defined by (\ref{Di}). The  vector $\textbf{W}_n $ is the similar of the vector  $\bfz_n$ when  $p$ was  fixed.  
Then, the process $G_n$ can be written:
\begin{equation}
\label{eq10}
G_n\pth{\ebo+C \sqrt{\frac{p}{n}} \eu} - G_n(\ebo)=\eE \cro{ G_n\pth{\ebo+C \sqrt{\frac{p}{n}} \eu} - G_n(\ebo) }+ \textbf{W}_n \eu +\sum^n_{i=1} [{\cal R}_i-\eE[{\cal R}_i]].
\end{equation}
First all, remark that 
\begin{equation}
\label{Lu}
\| \eu \|^2 \lambda_{\min} \big( \frac{1}{n} \eeX_i \eeX_i^t \big) \leq  \frac{1}{n} \sum^n_{i=1} \eu^t \eeX_i \eeX_i^t \eu \leq \| \eu \|^2 \lambda_{\max} \big( \frac{1}{n} \eeX_i \eeX_i^t \big). 
\end{equation}
Since the errors $(\varepsilon_i)_{1 \leqslant i \leqslant n}$  are  independent,  using also $|{\cal R}_i| < C \sqrt{\frac{p}{n}}\left|   \eeX^t_i \eu \right|\e1_{|\varepsilon_i| <  C \sqrt{\frac{p}{n}} \left| \eeX^t_i \eu \right| }$, together with assumption (A5), we obtain
\[
\eE \cro{\e1_{|\varepsilon_i| <  C \sqrt{\frac{p}{n}} \left| \eeX^t_i \eu \right| } } \leq C \sqrt{\frac{p}{n}} \| \eeX_i\|  \leq C \sqrt{\frac{p}{n}} \max_{1 \leqslant i \leqslant n} \| \eeX_i\| =o(1), 
\]
which imply, since  $(\varepsilon_i)_{1 \leqslant i \leqslant n}$ are i.i.d., 
\[
\eE \cro{ \sum^n_{i=1} [{\cal R}_i-\eE[{\cal R}_i]] }^2 = \sum^n_{i=1} \eE [{\cal R}_i-\eE[{\cal R}_i]]^2 \leq  \sum^n_{i=1} \eE [ {\cal R}^2_i]
\]
and by assumptions (A3), (A5) together with relation (\ref{Lu}), we have
\begin{equation}
\label{eq11}
\qquad \qquad  \qquad \qquad
\leq  C^2\frac{p}{n} \sum^n_{i=1} |\eeX^t_i \eu |^2 \eE \cro{\e1_{|\varepsilon_i| <  C \sqrt{\frac{p}{n}} \left| \eeX^t_i \eu \right| } }  = o \pth{ \frac{p}{n} \sum^n_{i=1} \eu^t \eeX_i \eeX_i^t \eu}=o(p).
\end{equation}
For the last relation we have used assumption (A4).\\
Let be the following random  variable   $U_n \equiv p^{-1/2} \sum^n_{i=1} [{\cal R}_i-\eE[{\cal R}_i]] $. Then, relation (\ref{eq11}) implies $\eE[U^2_n] =o(1)$. This, together with $\eE[U_n]=0$, imply, by the  Bienaymé-Tchebychev inequality, that  $U_n \overset{\PP} {\underset{n \rightarrow  \infty}{\longrightarrow}} 0 $. Thus $\sum^n_{i=1} [{\cal R}_i-\eE[{\cal R}_i]] =o_{\PP} (p^{1/2}) $. Then, relation (\ref{eq10}) becomes
\begin{equation}
\label{eq12}
G_n\pth{\ebo+C \sqrt{\frac{p}{n}} \eu} - G_n(\ebo)=\eE \cro{ G_n\pth{\ebo+C \sqrt{\frac{p}{n}} \eu} - G_n(\ebo) }+ \textbf{W}_n \eu + o_{\PP} (p^{1/2})
\end{equation}
which is equal to, using (\ref{eq9}):
\[
\qquad \qquad  \qquad  
=\cro{C^2 \frac{f(0)}{2} p \pth{ \frac{1}{n} \sum^n_{i=1} \eu^t \eeX_i \eeX_i^t \eu } +  C  \sqrt{p} \pth{ \sum^n_{i=1} \frac{{\cal D}_i \eeX^t_i}{\sqrt{n}}} \eu  } (1+o_{\PP}(1))+ o_{\PP} (p^{1/2}).
\]
Since $\eE[{\cal D}_i ]=0$, $Var[{\cal D}_i  \eeX_i^t \eu]= \tau (1- \tau) \eu^t \eeX_i \eeX_i^t \eu$ and $\| \eu \|=1$, then,  using (\ref{Lu}), we have that $n^{-1/2} \sum^n_{i=1} {{\cal D}_i \eeX^t_i \eu}  $ converges in law to a centred Gaussian distribution. Taking into account assumptions (A4) and (A6), for $B$ large enough, we obtain
\begin{equation}
\label{eq13}
G_n\pth{\ebo+B \sqrt{\frac{p}{n}} \eu} - G_n(\ebo)= B^2 f(0) p \pth{ \frac{1}{n} \sum^n_{i=1} \eu^t \eeX_i \eeX_i^t \eu } (1+o_{\PP}(1)) >0,
\end{equation}
for $n$ large enough. Thus, relation (\ref{eq8}) follows taking into account assumptions (A3) and (A4).
\hspace*{\fill}$\blacksquare$  \\

\noindent {\bf Proof of Theorem \ref{th_vconv}}.
We have the following inequality, with probability 1, for the quantile estimator $\widetilde{\eb}_{n;j}$:
\[
\min_{j \in {\cal A}} \| \eb^0_j\| \leq \max_{j \in {\cal A}} \|\widetilde{\eb}_{n;j}- \eb_j^0 \| +\min_{j \in {\cal A}} \| \widetilde{\eb}_{n;j}\|.
\]
By Lemma \ref{Lemma 2.1}, we have that $\max_{j \in {\cal A}} \|\widetilde{\eb}_{n;j}- \eb_j^0 \| =O_{\PP}((p/n)^{1/2}))=O_{\PP}( n^{(c-1)/2})$. On the other hand, we denoted in assumption (A8), $h_0=\min_{j \in {\cal A}} \| \eb^0_j\|$. Then, we have with probability one, $h_0n^{-\alpha} \leq O_{\PP}( n^{(c-1)/2 - \alpha})+n^{- \alpha} \min_{j \in {\cal A}} \| \widetilde{\eb}_{n;j}\|  $ and taking into account assumption (A8): $M \leq o_\PP (1)+n^{- \alpha}  \min_{j \in {\cal A}} \| \widetilde{\eb}_{n;j}\|$. Hence:
\begin{equation}
\label{L42}
\lim_{n \rightarrow \infty}\PP \bigg[ \min_{j \in {\cal A}} \|  \widetilde{\eb}_{n;j} \|>\frac{M n^\alpha}{2} \bigg] = 1.
\end{equation}
We consider the $r$-vector $\eu$ such that $\| \eu\|=1$. Similarly to  Theorem 2.1 of \cite{Zhang-Xiang-15}, we have, for a constant $B$,
\[
  \sum^p_{j=1} \|\widetilde{\eb}_{n;j} \|^{- \gamma} \cro{\|\eb^0_j+B\sqrt{\frac{p}{n}} \eu_j \|- \|\eb^0_j\|} \geq   \sum^{p_0}_{j=1} \|\widetilde{\eb}_{n;j} \|^{- \gamma} \cro{\|\eb^0_j+B\sqrt{\frac{p}{n}} \eu_j \|- \|\eb^0_j\|}
\geq - B\sqrt{\frac{p}{n}}   \sum^{p_0}_{j=1} \|\widetilde{\eb}_{n;j} \|^{- \gamma} \| \eu_j\|  .
\]
Since $\|\eu\|=1$, taking into account relation (\ref{L42}), applying the Cauchy-Schwarz inequality, we obtain
\[
- \sqrt{\frac{p}{n}} \lambda_n \sum^{p_0}_{j=1} \|\widetilde{\eb}_{n;j} \|^{- \gamma} \| \eu_j\| \geq - \sqrt{\frac{p}{n}} \lambda_n \big(\sum^{p_0}_{j=1} \|\widetilde{\eb}_{n;j} \|^{- 2\gamma} \big)^{1/2} \| \eu\| \geq  - \sqrt{\frac{p}{n}} \frac{\lambda_n \sqrt{p_0}}{\big( \min_{j \in {\cal A}} \|\widetilde{\eb}_{n;j} \| \big)^\gamma}  \geq - \sqrt{\frac{p}{n}} \frac{\lambda_n \sqrt{p_0}}{ \big( \frac{M n^\alpha}{2} \big)^\gamma}
\]
and by assumption $\lambda_n n^{(1+c)/{2} - \alpha \gamma}  {\underset{n \rightarrow \infty}{\longrightarrow}} 0$, we obtain:
\begin{equation}
\label{eq14}
\qquad \qquad  \qquad \qquad
= O_{\PP}\pth{\frac{p}{n} } .
\end{equation}
The Theorem is proved if we have the similar of inequality (\ref{eq8}) for
\[
 \frac{1}{n} G_n\pth{ \ebo+ B \sqrt{\frac{p}{n}} \eu} - \frac{1}{n} G_n(\ebo) +\sum^p_{j=1} \lambda_n \widehat{\omega}_{n;j} \cro{\|\eb^0_j+\sqrt{\frac{p}{n}} B \eu_j \|- \|\eb^0_j\| }
\equiv Q_n \pth{ \ebo+ B \sqrt{\frac{p}{n}} \eu} -Q_n(\ebo).
\]
We show that for all   $ \epsilon >0$ there exists  $  B_\epsilon$ large enough, such that, for any $n$ large enough
\begin{equation}
\label{eq8bis}
\PP \cro{ \inf_{\| \eu \|=1}Q_n \pth{ \ebo+ B_\epsilon \sqrt{\frac{p}{n}} \eu}  > Q_n(\ebo) } > 1-\epsilon .
\end{equation}
By the  definition of $Q_n$, we have for all constant $B>0$, that
\[
Q_n \pth{ \ebo+ B \sqrt{\frac{p}{n}} \eu} -Q_n(\ebo)  
= \frac{1}{n} G_n\pth{ \ebo+ B \sqrt{\frac{p}{n}} \eu} - \frac{1}{n} G_n(\ebo) +\sum^p_{j=1} \lambda_n \widehat{\omega}_{n,j} \cro{\|\eb^0_j+\sqrt{\frac{p}{n}} B \eu_j \|- \|\eb^0_j\| }
\]
and using  relations (\ref{eq13}) and (\ref{eq14})
\[
\qquad \qquad  \qquad \qquad
> B^2 f(0) \frac{p}{n} \pth{\frac{1}{n}  \sum^n_{i=1} \eu^t \eeX_i \eeX_i^t \eu } (1+o_{\PP}(1)) - B O_{\PP}\pth{ \frac{p}{n}}  .
\]
Relation (\ref{eq8bis}) follows from the last relation,   for $n$ and $B$ large enough and using assumption (A4).
\hspace*{\fill}$\blacksquare$ \\

\noindent {\bf Proof of Theorem \ref{Theorem 2SPL}}. 
\textit{(i)} 
By  Theorem \ref{th_vconv}, we have that $\widehat{\eb}^*_n$ belongs,  with a probability converging to one, to the set: ${\cal V}_p(\ebo) \equiv \acc{ \eb; \| \eb-\ebo\| \leq B \sqrt{\frac{p}{n}}}$, with $B>0$ large enough as in relation (\ref{eq8bis}).\\
\textit{For $p >p_0$}, we show that for all $\eb=(\eb_{\cal A}, \eb_{{\cal A}^c}) \in {\cal V}_p(\ebo)$ such that $\| \eb_{\cal A} - \eb^0_{\cal A} \| = O\pth{ \sqrt{\frac{p}{n}}}$ and for all constant $C \in (0, B)$, we have
\begin{equation}
\label{eq15}
Q_n(\eb_{\cal A}, \textbf{0}_{p-p_0})=\min_{\|\eb_{{\cal A}^c}\| \leq C \sqrt{\frac{p}{n}}} Q_n(\eb_{\cal A}, \eb_{{\cal A}^c}),
\end{equation}
with a probability tending to one, as sample size $n \rightarrow \infty$.\\
Let us consider the parameter set ${\cal W}_n \equiv \acc{\eb \in {\cal V}_p(\ebo) ; \|\eb_{{\cal A}^c}\| >0}$.\\ 
We show that $\PP[\widehat{\eb}^*_n \in {\cal W}_n] \rightarrow 0$, as $n \rightarrow \infty$. 
For this, we firstly consider  two parameter  vectors $\eb=(\eb_{\cal A}, \eb_{{\cal A}^c}) \in {\cal W}_n$ and $ {\eb}^{(1)}=({\eb}_{\cal A}^{(1)}, {\eb}_{{\cal A}^c}^{(1)}) \in {\cal V}_p(\ebo)$,   such that ${\eb}^{(1)}_{\cal A}=\eb_{\cal A}$ and $ {\eb}_{{\cal A}^c}^{(1)}=\textbf{0}_{d-d^0}$.\\
Let us take the difference of the objective random process for the two parameter vectors. We denote this difference $D_n (\eb, \eb^{(1)})$:
\begin{equation}
\begin{array}{cl}
D_n (\eb, \eb^{(1)}) & \equiv Q_n(\eb) - Q_n({\eb}^{(1)})
   \\
& = n^{-1} \sum^n_{i=1} \cro{\rho_\tau(Y_i - \eeX^t_i \eb)-\rho_\tau(Y_i - \eeX^t_i {\eb}^{(1)})} +\sum^p_{j=p_0+1} \lambda_n \widehat{\omega}_{n;j} \| \eb_j\|.
\end{array} 
\label{DD}
\end{equation}
From \cite{Knight-98}, we have the following identity, for any $x, y \in \R$:  
\[
\rho_\tau(x-y)- \rho_\tau(x)=y(\e1_{x \leq 0} - \tau)+\int^y_0 (\e1_{x \leq t} -\e1_{x \leq 0})dt.
\]
Using this  relation for the first sum of (\ref{DD}), we obtain:
\begin{equation}
\begin{array}{c}
\displaystyle{ \frac{1}{ n} \sum^n_{i=1} \big[\rho_\tau(Y_i- \eeX^t_i \eb )-\rho_\tau(Y_i - \eeX^t_i {\eb}^{(1)} )\big]= \frac{1}{ n} (\eb  -{\eb}^{(1)} )^t \sum^n_{i=1} \eeX_i  [\e1_{Y_i - \eeX^t_i {\eb}^{(1)} \leq 0} - \tau]}  \\
\displaystyle  {+\frac{1}{n} \sum^n_{i=1} \int^{\eeX_i^t (\eb  - {\eb}^{(1)} )}_0 [\e1_{Y_i - \eeX^t_i {\eb}^{(1)}  \leq t} -\e1_{Y_i - \eeX^t_i {\eb}^{(1)} \leq 0}] dt} 
\equiv T_{1n}+T_{2n} .
\end{array} 
\label{eqq0}
\end{equation}
For $T_{1n}$, since the density $f$ is bounded in a neighbourhood of  0,  we have, by assumption (A5), that:
\[
\eE[T_{1n}]= (\eb - {\eb}^{(1)} )^t \frac{1}{ n} \sum^n_{i=1} \eeX_i [F(\eeX^t_i ({\eb}^{(1)}  -\eb^0 ))-F(0)] = (\eb - {\eb}^{(1)} )^t \frac{f(0)}{ n}  \sum^n_{i=1}\eeX_i  \eeX^t_i (\eb^0  - {\eb}^{(1)})  (1+o(1)).   
\]
Then, 
$
 |\eE[T_{1n}]| \leq \|\eb - {\eb}^{(1)}  \|  \cdot \big\| n^{-1}\sum^n_{i=1}\eeX_i  \eeX^t_i  \big\|\cdot \| \eb^0  - {\eb}^{(1)}\| f(0) (1+o(1)) $.
 Since the matrix  $ n^{-1}\sum^n_{i=1}\eeX_i  \eeX^t_i $ is Hermitian, we have, taking into account assumption (A4), that $n^{-1} \left\| \sum^n_{i=1}\eeX_i  \eeX^t_i  \right\| $ $= \lambda_{\max}(n^{-1}\sum^n_{i=1}\eeX_i  \eeX^t_i )$ $ \leq M_0$. Hence,   we have 
$
|\eE[T_{1n}]| \leq M_0 f(0) \|\eb - {\eb}^{(1)}  \|  \cdot \| \eb^0  - {\eb}^{(1)}\| $.
Therefore
$
\eE[T_{1n}]=O( \|\eb - {\eb}^{(1)}  \|^2 )$. 
By calculations analogous to $\eE[T_{1n}]$, using independence of $(\varepsilon_i)_{1 \leqslant i \leqslant n}$, we have that 
$\eE[T^2_{1n}] = C n^{-1} {\|\eb - {\eb}^{(1)} \|^3} \rightarrow 0$,  for $n \rightarrow \infty$. Since $Var[T_{1n}] \leq \eE[T^2_{1n}]$, using the Bienaymé-Tchebychev inequality, we have:
\begin{equation}
\label{eqq1}
T_{1n}= C \|\eb - \widetilde \eb \|^2 (1+o_{\PP}(1)).
\end{equation}
Consider now   $T_{2n}$ of relation (\ref{eqq0}), which can be written as:
\[
T_{2n}= n^{-1} \sum^n_{i=1} \int^{\eeX^t_i (\eb  - {\eb}^{(1)} )}_0 [\e1_{\varepsilon_i \leq t- \eeX_i^t (\eb^0  - {\eb}^{(1)} )} - \e1_{\varepsilon_i \leq - \eeX_i^t (\eb^0  - {\eb}^{(1)})} ] dt .\]
Taking into account that $\eb  \in {\cal V}_{p}(\eb^0_n) $, together with assumptions (A3),  (A5), we get 
$
\eE[T_{2n}]  =  n^{-1} \sum^n_{i=1} \int^{\eeX^t_i (\eb  - {\eb}^{(1)} )}_0[F(t- \eeX_i^t (\eb^0  - {\eb}^{(1)} )) - F(- \eeX_i^t (\eb^0  - {\eb}^{(1)}))] dt  
  =  n^{-1} \sum^n_{i=1} \int^{\eeX^t_i (\eb  - {\eb}^{(1)} )}_0 [t  f(\eeX_i^t ({\eb}^{(1)} $ $ - \eb^0 )) + o(t) ]dt $.   \\ 
 By the proof of  Lemma \ref{Lemma 2.1},  assumptions (A3),   (A5), we obtain that $f(\eeX_i^t ({\eb}^{(1)}  - \eb^0 ))$ is bounded by a constant $C >0$. 
 Thus, as for $T_{1n}$, using assumption (A4) and the fact that 
 $
n^{-1} \sum^n_{i=1} \| \eeX_i\|^2 -tr\big( n^{-1}  \sum^n_{i=1} \eeX_i \eeX_i^t\big)  {\underset{n \rightarrow \infty}{\longrightarrow}}  0, 
 $
 we have
 $
 |\eE[T_{2n}]| \leq 
  {C}n^{-1} \sum^n_{i=1} \| \eeX_i \|^2   \| \eb  -{\eb}^{(1)}  \| \cdot \| {\eb}^{(1)}  - \eb^0 \|   +o\big( n^{-1} \sum^n_{i=1}$ $\eeX_i^t (\eb -{\eb}^{(1)} )\big)  
 = C \| \eb  -{\eb}^{(1)}  \|^2  $.
We show similarly that $\displaystyle{\eE[T_{2n}^2] = C n^{-1} \| \eb  - {\eb}^{(1)} \|^3}$. 
Then, by the  Bienaymé-Tchebychev inequality, we get:
\begin{equation}
\label{eqq2}
T_{2n}= C \| \eb - {\eb}^{(1)} \|^2 (1+o_{\PP}(1)).
\end{equation}
Hence, by relations  (\ref{eqq1}), (\ref{eqq2}), we obtain
\begin{equation}
\label{TT}
T_{1n}+T_{2n} = C \| \eb  - {\eb}^{(1)} \|^2 (1+o_{\PP}(1)).
\end{equation}
This last  relation together with relations   (\ref{eqq0}), (\ref{TT}), give for relation  (\ref{DD}):
\[
D_n (\eb, {\eb}^{(1)})=C \| \eb-{\eb}^{(1)}\|^2 (1+o_{\PP}(1))+\sum^p_{j=p_0+1} \lambda_n \widehat{\omega}_{n;j} \| \eb_j\|.
\]
On the other hand, by Lemma \ref{Lemma 2.1}, we have
\[
\widehat{\omega}_{n;j}=\frac{1}{\| \widetilde{\eb}_{n;j} \|^\gamma} =\frac{1}{\| \widetilde{\eb}_{n;j} -{\eb}^0_{j} \|^\gamma} =O_{\PP} \pth{\pth{\frac{p}{n}}^{- \gamma /2}}
\]
and moreover  for all $j \geq p_0+1$, since $\eb \in {\cal W}_n \subseteq {\cal V}_p(\ebo)$ we have $0< \| \eb_j\| =O\pth{ \sqrt{\frac{p}{n}}}$. Then 
\[
\sum^p_{j=p_0+1} \lambda_n \widehat{\omega}_{n;j} \| \eb_j\| = \sum^p_{j=p_0+1} \lambda_n O_{\PP} \pth{ \pth{\frac{p}{n}}^{(1-\gamma)/2}}.
\]
Thus
\[
\frac{D_n(\eb, {\eb}^{(1)})}{\| \eb -{\eb}^{(1)} \|} \geq C \| \eb -{\eb}^{(1)} \| (1+o_{\PP}(1))+   \sum^p_{j=p_0+1} \lambda_n O_{\PP} \pth{ \pth{\frac{p}{n}}^{-\gamma/2}}.
\]
We have that $\| \eb - {\eb}^{(1)} \|= O\pth{ \pth{\frac{p}{n}}^{1/2} }$. Since ${p}/{n}=O(n^{c-1})$, under the assumption that $\lambda_n n^{(1-c){(1+\gamma)}/{2}} \rightarrow \infty $, as $n \rightarrow\infty$, and $p>p_0$, we have  
\begin{equation}
\label{dd1}
\frac{D_n(\eb, {\eb}^{(1)})}{\| \eb - {\eb}^{(1)} \|} \geq \lambda_n O_{\PP} \left( n^{{\gamma (1-c)}/{2}} \right).
\end{equation}
To finish the proof of relation (\ref{eq15}),  consider now  other two parameter vectors: $\ebo$ the true value and ${\eb}^{(1)}$ a parameter such that ${\eb}^{(1)} \equiv ({\eb}^{(1)}_{\cal A}, {\eb}^{(1)}_{{\cal A}^c}  )$, ${\eb}^{(1)}_{\cal A}={\eb}_{\cal A}$, $ {\eb}^{(1)}_{{\cal A}^c} = \eb^0_{{\cal A}^c} = \textbf{0}_{d-d^0}$. We obtain as for (\ref{TT}) that:
\[
D_n(\ebo, {\eb}^{(1)}) = n^{-1} \sum^n_{i=1} \cro{ \rho_\tau(Y_i - \eeX^t_i \ebo )- \rho_\tau(Y_i - \eeX^t_i {\eb}^{(1)} ) }=C \| (\ebo - {\eb}^{(1)})_{{\cal A}} \|^2 (1+o_{\PP}(1)).
\]
Then
\begin{equation}
\label{dd2}
\frac{D_n(\ebo, {\eb}^{(1)})}{\| \ebo -{\eb}^{(1)} \|} = C \| (\ebo - {\eb}^{(1)})_{{\cal A}} \| (1+o_{\PP}(1)) = O_{\PP}\pth{ \pth{\frac{p}{n}}^{1/2} } = O_{\PP} \pth{ n^{(c-1)/2}}.
\end{equation}
Since $\lambda_n n^{(1-c){(1+\gamma)}/{2}} \rightarrow \infty$ as $n \rightarrow \infty$, we have that  (\ref{dd1}) is much bigger than (\ref{dd2}), for $n$ large enough. Then,  relation (\ref{eq15}) follows. \\
To finish the proof of \textit{(i)}, we will show that
\begin{equation}
\label{Pmb}
\lim_{n \rightarrow \infty}\PP \bigg[   \min_{j \in {\cal A}} \| \widehat{\eb}^*_{n;j} \| > 0  \bigg] = 1.
\end{equation}
With probability 1, we have that:
$
\min_{j \in {\cal A}} \| \widehat{\eb}^*_{n;j} \| \geq \min_{j \in {\cal A}} \|  {\eb}^0_{j} \| -  \max_{j \in {\cal A}} \|\widehat{\eb}^*_{n;j} -  {\eb}^0_{j}  \|$. 
By Theorem \ref{th_vconv} and assumption (A8),  we have,  
\[
\lim_{n \rightarrow \infty}\PP \bigg[n^{- \alpha}  \min_{j \in {\cal A}} \| \widehat{\eb}^*_{n;j} \| \geq \frac{M}{2}  \bigg] = 1, 
\]
which implies
\[
\lim_{n \rightarrow \infty}\PP \bigg[n^{- \alpha}  \min_{j \in {\cal A}} \| \widehat{\eb}^*_{n;j} \| > 0  \bigg] = 1,
\]
from which relation (\ref{Pmb}) follows. Relations (\ref{Pmb}) and (\ref{eq15}) imply \textit{(i)}.\\
\textit{(ii)} 
Taking into account claim \textit{(i)} and assumption (A7), the estimator $\widehat{\eb}^*_n$ can be written with a probability converging to 1 as $\widehat{\eb}^*_n = \ebo+\sqrt{\frac{p}{n}} \ed$, with  $\ed=(\ed_{\cal A}, \ed_{{\cal A}^c})$, $\ed_{{\cal A}^c}=\textbf{0}_{r-r^0}$, $\| \ed_{\cal A}\| \leq C $. Then
\begin{equation}
\label{eq16}
Q_n \pth{\ebo+ \sqrt{\frac{p}{n}} \ed}  - Q_n(\ebo)=\frac{1}{n} \sum^n_{i=1} \cro{\rho_\tau\pth{Y_i-\eeX^t_i(\ebo+\sqrt{\frac{p}{n}} \ed)}-\rho_\tau(\varepsilon_i) }+{\cal P},
\end{equation}
with ${\cal P} \equiv \sum^{p_0}_{j=1} \lambda_n \widehat{\omega}_{n;j} \cro{ \| \eb_j\| - \| \eb^0_j\|}$.\\
For all $j \in \{ 1, \cdots , p_0\}$ we have $\left| \| \widehat{\eb}^*_{n;j} \| - \| \eb^0_j \| \right|=C \sqrt{\frac{p}{n}} \ed^t_j \eb^0_j $, with a probability converging to  1. On the other hand, using assumption (A8) and Lemma \ref{Lemma 2.1} we have for the quantile estimators
\[
\| \widetilde {\eb}_{n;j}\| =\| \ebo_j\|+ O_{\PP} \pth{\sqrt{\frac{p}{n}}} \geq h_0 +O_{\PP} \pth {n^{- (1-c)/{2}}}= O(n^\alpha) +O_{\PP} \pth {n^{- (1-c)/{2}}}= O(n^\alpha).
\]
Thus, we have  $\widehat{\omega}_{n;j} \leq O_{\PP}(n^{-\alpha \gamma})$, for all  $j \in {\cal A}$.
Then,  for the  the second term on the right hand side of (\ref{eq16}), we have with a probability converging to one: 
\begin{eqnarray}
|{\cal P}| &= & \sum^{p_0}_{j=1 }\lambda_n \widehat{\omega}_{n;j} \left| \| \eb^0_j +\sqrt{\frac{p}{n}} \ed_j \| - \| \eb^0_j \|  \right|  
 \leq   C n^{-\alpha \gamma} \sum^{p_0}_{j=1 } \lambda_n \sqrt{\frac{p}{n}} \left|\ed_j^t \eb^0_j \right| \nonumber \\
 & \leq & C n^{-\alpha \gamma} \sqrt{\frac{p}{n}} \lambda_n \sum^{p_0}_{j=1 } \| \ed_j \| \cdot \| \eb^0_j\|   
 \leq  C \sqrt{\frac{p}{n}}  r^0 \lambda_n n^{-\alpha \gamma} = C r^0 \lambda_n n^{(c-1)/{2} - \alpha \gamma}.
 \label{eq17}
\end{eqnarray}
We study now the first term of the right hand side of  (\ref{eq16}), which can be written as:
\[
\frac{1}{n} \sum^n_{i=1} \cro{\rho_\tau \pth{ Y_i - \eeX^t_i \pth{ \ebo + \sqrt{\frac{p}{n}} \ed}}-  \rho_\tau(\varepsilon_i) }  \qquad \qquad  \qquad \qquad \qquad \qquad
\]
\[
= \frac{1}{n} \sqrt{\frac{p}{n}} \sum^n_{i=1} \eeX^t_i \ed [\e1_{\varepsilon \leq 0} - \tau] +\frac{1}{n}\sum^n_{i=1} \int^{\sqrt{\frac{p}{n}} \eeX^t_i \ed}_ 0 [ \e1_{\varepsilon_i \leq t} - \e1_{\varepsilon_i \leq 0} ] dt \equiv J_1+J_2.
\]
Since $\eE[J_1]=0$,  using assumption (A5) and  the  Cauchy-Schwarz inequality, we have:
\[
Var[J_1] \leq \eE[J_1^2]= \frac{1}{n^2} \frac{p}{n} \tau(1-\tau)\sum^n_{i=1} (\eeX^t_i \ed )^2 \leq \frac{p}{n^3} \tau(1-\tau) \sum^n_{i=1}   \|\eeX_{i, {\cal A}} \|^2 \| \ed_{\cal A} \|^2 \leq C \frac{p}{n^3} n  \max_{1 \leq i  \leq n} \|\eeX_{i, {\cal A}} \|^2   {\underset{n \rightarrow \infty}{\longrightarrow}}  0.
 \]
Using assumption (A2), we have for the expectation of $J_2$:
\begin{equation} 
\label{EJ2}
\eE[J_2]=\frac{1}{n} \int^{\sqrt{\frac{p}{n}} \eeX^t_i \ed}_0 (t f(0)+o(t^2))dt =  \frac{f(0)}{2} \frac{p}{n} \ed^t \frac{1}{n}\sum^n_{i=1} \eeX_i \eeX_i^t \ed (1+o(1)).
\end{equation}
Using assumption (A4) and relation (\ref{Lu}), we have that
\[
f(0)p n^{-1}  \|\ed\|^2_2   \lambda_{\min}\big( n^{-1} \sum^n_{i=1}\eeX_{i  } \eeX^t_{i  }   \big) \leq \eE[J_2] \leq f(0)p n^{-1}  \|\ed\|^2_2   \lambda_{\max}\big(n^{-1} \sum^n_{i=1}\eeX_{i } \eeX^t_{i  }   \big).\] 
Taking into account the fact that $\|\ed\|^2_2=\|\ed_{{\cal A}}\|^2_2 \leq C  $, we have:
\begin{equation}
\label{eq18}
\eE[J_2] = C f(0) \frac{p}{n}.
\end{equation}
Similarly we obtain $Var[J_2]=O\pth{ \frac{1}{n} \pth{\frac{p}{n}}^{3/2}  } \rightarrow 0$, as $n \rightarrow \infty$. \\
But  $\eE[J_2]=O\pth{n^{c-1}}$ and by assumption (A7), together with  relation (\ref{eq17}),  we have $|{\cal P}| \leq O_{\PP} \pth{\lambda_n n^c n^{(c-1)/{2} - \alpha \gamma}} = O_{\PP} \pth{\lambda_n   n^{(3c-1)/{2} - \alpha \gamma}} $. Then $\frac{|{\cal P}| }{\eE[J_2]} = O_{\PP} \pth{\lambda_n   n^{(c+1)/{2} - \alpha \gamma}} \rightarrow 0$ and thus  $\eE[J_2] \gg |{\cal P}|$. \\
Hence, minimizing (\ref{eq16}) amounts to  minimizing $J_1+J_2$, with respect to $ \sqrt{\frac{p}{n}} \ed$. Using relation (\ref{EJ2}), we obtain:
\[
\frac{1}{n} \sum^n_{i=1} [\rho_\tau(Y_i -\eeX_i^t( \eb^0_n + \sqrt{\frac{p}{n}} \ed)) - \rho_\tau(\varepsilon_i)] = \frac {\sqrt{\frac{p}{n}}}{ n}  \sum^n_{i=1}\eeX^t_{i,{{\cal A}}} \ed_{{\cal A}} [\e1_{\varepsilon_i <0} - \tau]+  \frac{f(0)}{2} \frac{p}{n} \ed^t_{{\cal A}} \eU_{n,{\cal A}} \ed_{{\cal A}}(1+o_{\PP}(1)).
\]
The minimizer of the right hand side of the last equation is:
\begin{equation}
\label{ad}
 \sqrt{\frac{p}{n}} \ed_{{\cal A}}=-\frac{1}{n} \frac{1}{f(0)} \eU^{-1}_{n,{\cal A}} \big( \sum^n_{i=1} \eeX_{i,{ {\cal A}}} (\e1_{\varepsilon_i \leq 0}-\tau)\big).
\end{equation}
For studying (\ref{ad}), let us consider the following   independent random variable sequence, for $i = 1, \cdots, n$, 
\[
W_i \equiv (f(0))^{-1} \eu^t \eU^{-1}_{n,{\cal A}}\eeX_{i {\cal A}} (\e1_{\varepsilon_i \leq 0}-\tau),\] 
with $\eu $ a  vector of  dimension $r^0$,  such that  $\|\eu\|_2=1$. We have that  $\eE[W_i]=0$ and
$
\sum^n_{i=1}Var[W_i]=n \tau(1-\tau)(f(0))^{-2}\eu^t \eU^{-1}_{n,{\cal A}} \eu$.
Then, by CLT for  independent random variable sequences $(W_i)_{1 \leqslant i \leqslant n}$, we have
\begin{equation}
\label{cv}
\sqrt{n} f(0) \frac{\eu^t (\widehat \eb_{{\cal A}}- \eb^0_{{\cal A}})}{\sqrt{\tau(1- \tau) (\eu^t \eU^{-1}_{n,{\cal A}}  \eu)}}  \overset{\cal L} {\underset{n \rightarrow \infty}{\longrightarrow}} {\cal N}(0,1).
\end{equation}
Claim (ii) results by taking into account of the fact that  $\widehat \eb_{{\cal A}}- \eb^0_{{\cal A}}= \sqrt{\frac{p}{n}} \ed_{{\cal A}}$, together with relations  (\ref{ad}), (\ref{cv}).

\hspace*{\fill}$\blacksquare$ \\


\begin{thebibliography}{3}
\bibitem[Ciuperca(2015)]{Ciuperca-15}
Ciuperca, G., (2015).  
\newblock Model selection in high-dimensional quantile regression   with  seamless $L_0$ penalty. 
\newblock{\it  Statistics and Probability Letters},   {\bf 107},  313-323. 
\bibitem[Ciuperca(2016)]{Ciuperca-15b}
Ciuperca, G., (2016). 
\newblock Adaptive LASSO model selection in a multiphase quantile regression. 
\newblock{\it  Statistics}, DOI: 10.1080/02331888.2016.1151427.
\bibitem[Guo et al.(2015)]{Guo-Zhang-Wang-Wu-15}
Guo, X., Zhang, H., Wang, Y., Wu, J.L., (2015). 
\newblock Model selection and estimation in high dimensional regression models with group SCAD. 
\newblock{\it  Statistics and Probability Letters},   {\bf 103},  86-92. 
\bibitem[Huang et al.(2012)]{Huang-Breheny-Ma-12}
Huang, J., Breheny, P., Ma, S., (2012). 
\newblock A selective review of group selection in high-dimensional models.
\newblock{\it  Statistical Science},  {\bf 27}, No. 4,   481-499.
\bibitem[Knight (1998)]{Knight-98}
Knight,  K., (1998). 
\newblock Limiting distributions for L1 regression estimators under general conditions.
\newblock {\it  The Annals of Statistics},  {\it 26},  755-770.
\bibitem[Koenker(2005)]{Koenker-05}
Koenker, R., (2005).
\newblock{\it Quantile Regression}. 
\newblock{Cambridge University Press}.
\bibitem[Nardi and Rinaldo(2008)]{Nardi-Rinaldo-08}
Nardi, Y., Rinaldo, A., (2008). 
\newblock On the asymptotic properties of the group lasso estimator for linear models.
\newblock{\it  Electronic Journal of Statistics}, {\bf 2},  605-633.
\bibitem[Xu and Ghosh(2015)]{Xu-Ghosh-15}
Xu, X., Ghosh, M., (2015). 
\newblock Bayesian variable selection and estimation for group Lasso.
\newblock{\it Bayesian Analysis}, {\bf 10}(4),  909-936.
\bibitem[Yuan and Lin(2008)]{Yuan-Lin-06}
Yuan, M., Lin, Y., (2006).  
\newblock Model selection and estimation in regression with grouped variables.
\newblock{\it J. R. Statist. Soc. B}, {\bf 68}(1),  49-67.
\bibitem[Zhang and Xiang(2015)]{Zhang-Xiang-15}
Zhang, C., Xiang, Y., (2015). 
\newblock On the oracle property of adaptive group LASSO in high-dimensional linear models.  
\newblock{\it Statistical Papers}, \textbf{57(1)}, 249-265.
\bibitem[Zou and Yuan(2008)]{Zou-Yuan-08}
Zou, H., Yuan, M., (2008). 
\newblock Composite quantile regression and the oracle model selection theory.
\newblock{\it The Annals of Statistics}, {\bf 36}(3),  1108-1126.
\bibitem[Zou and Zhang(2009)]{Zou-Zhang-09}
Zou, H., Zhang, H.H., (2009). 
\newblock On the adaptive elastic-net with a diverging number of parameters.
\newblock{\it The Annals of Statistics}, {\bf 37}(4),  1733-1751.
\bibitem[Wang and Leng(2008)]{Wang-Leng-08}
Wang, H., Leng, C., (2008). 
\newblock A note on adaptive group lasso.
\newblock{\it Computational Statistics and Data Analysis}, {\bf 52},  5277-5286.
\bibitem[Wang et al.(2015)]{Wang-You-Lian-15}
Wang, L., You, Y., Lian, H., (2015). 
\newblock Convergence and sparsity of Lasso and group Lasso in high-dimensional generalized linear models.
\newblock{\it Statistical Papers}, {\bf 56}, No. 3,   819-828.
\bibitem[Wei and Huang(2010)]{Wei-Huang-10}
Wei, F., Huang, J., (2010). 
\newblock Consistent group selection in high-dimensional linear model.
\newblock{\it Bernoulli}, {\bf 16}(4),  1369-1384.
\bibitem[Wu and Liu (2009)]{Wu-Liu-09}
Wu, Y., Liu, Y., (2009). 
\newblock Variable selection in quantile regression.
\newblock {\it Statistica Sinica}, {\bf 19},  801-817.
\end{thebibliography}
\end{document}